\documentclass[12pt]{amsart}
\usepackage[margin=0.8in]{geometry}
\newcommand{\R}{\mathbb{R}}

\newcommand{\N}{\mathbb{N}}

\usepackage{eucal}

\newcommand{\NN}{{{\mathcal{N}}}}

\renewcommand{\epsilon}{\varepsilon}

\newcommand{\e}{\varepsilon}

\usepackage[english]{babel}
\usepackage[alphabetic]{amsrefs}
\usepackage{amsmath,amssymb,amsfonts,amsthm,enumerate}
\usepackage{url}
\usepackage{graphicx,epstopdf,color}

\numberwithin{equation}{section}

\newtheorem{theorem}{Theorem}[section]
\newtheorem{corollary}[theorem]{Corollary}
\newtheorem{lemma}[theorem]{Lemma}

\newtheorem{remark}[theorem]{Remark}

\renewcommand{\leq}{\leqslant}
\renewcommand{\le}{\leqslant}
\renewcommand{\geq}{\geqslant}
\renewcommand{\ge}{\geqslant}

\begin{document}

\title[Fractional time-diffusion]{Decay estimates for evolutionary equations \\ with fractional time-diffusion}

\author[Serena Dipierro, Enrico Valdinoci and Vincenzo Vespri]{
Serena Dipierro${}^{(1)}$
\and
Enrico Valdinoci${}^{(1,2)}$
\and Vincenzo Vespri${}^{(3)}$
}
\subjclass[2010]{26A33, 34A08,	35K90	47J35, 58D25.}
\keywords{Fractional diffusion, parabolic equations, decay of solutions in time with respect to Lebesgue norms.}
\thanks{This work has been carried out during a very pleasant visit of the third author
at the School of Mathematics and Statistics at the University of Melbourne.
Supported by the Australian Research Council Discovery Project grant ``N.E.W. {\it Nonlocal
Equations at Work}''
and by the G.N.A.M.P.A. project ``{\it Nonlocal and degenerate
problems in the Euclidean space}''. The authors are members of~G.N.A.M.P.A.--I.N.d.A.M.}

\maketitle

{\scriptsize \scriptsize \begin{center} (1) -- 
Department of Mathematics and Statistics  \\
University of Western Australia\\
35 Stirling Hwy, Crawley, WA-6009 Perth
(Australia).\\
\end{center}\begin{center}(2) -- 
Dipartimento di Matematica
``Federigo Enriques''\\
Universit\`a
degli studi di Milano\\
Via Saldini 50, I-20133 Milano (Italy).\\
\end{center}
\scriptsize \begin{center} (3) -- 
Dipartimento di Matematica e Informatica ``Ulisse Dini''\\
Universit\`a
degli studi di Firenze\\
Viale Morgagni 67/a, I-50134 Firenze (Italy).\\ 
\end{center} 
\bigskip

\begin{center}
E-mail addresses: 
serena.dipierro@uwa.edu.au,
enrico.valdinoci@uwa.edu.au,
vincenzo.vespri@unifi.it
\end{center}
}

\begin{abstract}
We consider an evolution equation whose time-diffusion is of fractional type
and we provide decay estimates in time for the $L^s$-norm of the solutions
in a bounded domain.
The spatial operator that we take into account is very general and comprises
classical local and nonlocal diffusion equations.
\end{abstract}

\section{Introduction}

\subsection{General time-fractional diffusion equations}
The goal of this paper is to consider evolutionary equations
with nonlocal time-diffusion of fractional type, which is modeled by
an integro-differential operator. The space-diffusion
that we take into account can be both local and nonlocal, and in fact our
approach aims at general energy estimates in an abstract framework which
will in turn provide asymptotic decay
estimates in a series of concrete cases, including nonlocal nonlinear operators,
nonlocal porous medium equations and possibly nonlocal mean curvature operators.

More specifically, we consider equations of the form
\begin{equation}\label{0-2iAM-x12:A} \partial_t^{\alpha} u+\NN[u]=0,\end{equation}
with~$\alpha\in(0,1)$. In this setting,
the solution~$u$ is a function~$u=u(x,t)$,
with~$x$ lying in a nice Euclidean domain, $t>0$,
and Dirichlet boundary data. The variable~$x$ will be referred to as ``space'',
and, in the examples that we take into account, the operator~$\NN$
possesses some kind of ``elliptic'' features, which make~\eqref{0-2iAM-x12:A}
a sort of ``diffusive'', or ``parabolic'', equation.
In this spirit,
the variable~$t$, which will be referred to as ``time'',
appears in~\eqref{0-2iAM-x12:A} with a fractional derivative of
order~$\alpha\in(0,1)$, and we thereby consider~\eqref{0-2iAM-x12:A} as
a fractional time-diffusion.
In the examples that we take into account,
the diffusion modeled on the operator~$\NN$
can be either ``classical'' (i.e. involving derivatives of integer order,
up to order two), or ``anomalous'' (since it can involve fractional derivatives
as well, in which case we refer to it as a fractional space-diffusion).
\medskip

We also recall that integro-differential equations
are a classical topic in mathematical analysis, see e.g.~\cites{PIN, VOL}.
Fractional calculus also appears under different forms
in several real-world phenomena, see e.g.~\cites{ME, TRU, VU}.
In particular, time-fractional derivatives find applications
in the magneto-thermoelastic
heat conduction~\cite{MAG}, wave equations~\cites{CHEN, X9X},
hydrodynamics~\cite{MASS}, quantum physics~\cite{BH}, etc.
See also~\cites{LUC, TOP} for existence and uniqueness results and~\cites{VER0, AL2, AL1}
for related regularity results in the local and nonlocal spatial regime.
The recent literature has also widely considered time-fractional
diffusion coupled with $p$-Laplacian space-diffusion, see
e.g.~\cites{CHEN2, LI, LIU, ZH}
and the references therein.

In the framework of nonlocal equations, a deep and useful setting
is that provided by the Volterra integral equations, which
often offers a general context in which one develops existence, uniqueness,
regularity and asymptotic theories, see for
instance~\cites{RV0, RV1, RV2, RV22, RV3, RV4, RV5, RV55, RV6, RV7} and the references therein.
\medskip

The setting in which we work in this paper is the following. 
We consider the so-called Caputo derivative of 
order~$\alpha\in(0,1)$, defined as
\begin{equation} \label{K-DER}
\partial_t^{\alpha} u(t):=
\frac{d}{dt}
\int_0^t \frac{ u(\tau)-u(0)}{ (t-\tau)^\alpha }\,d\tau,
\end{equation}
up to a positive normalization constant, that we neglect (see e.g.~\cite{CAPU}).

Goal of this paper is to study
solutions
$$u=u(x,t):\R^n\times[0,+\infty)\to[0,+\infty)$$
of the initial value problem
\begin{equation}\label{IVP}
\left\{\begin{matrix}
\partial_t^{\alpha} u(x,t)+\NN[u](x,t)=0 & {\mbox{ for any }}x\in\Omega {\mbox{ and }}t>0, \\
u(x,t)=0 & {\mbox{ for any }}
x\in\R^n\setminus\Omega {\mbox{ and }}t\ge0,\\
u(x,0)=u_0(x) & {\mbox{ for any }}x\in\Omega.
\end{matrix}\right.\end{equation}
In our notation, $\Omega$ is a bounded subset of~$\R^n$ with smooth boundary
and~$\NN$ is a possibly nonlinear operator.
For concreteness, we suppose that the initial datum~$u_0$ does not vanish identically and lies
in~$L^q(\Omega)$ for any~$q\in[1,+\infty)$ (as a matter
of fact, weaker assumptions can be taken according to
suitable choices of the parameters). In any case,
from now on, the initial datum~$u_0$ will be always implicitly
supposed to be nonnegative, nontrivial and integrable at any power,
and the solution~$u$ to be nonnegative and smooth.\medskip

The main structural assumption that we take is that there exist~$s\in(1,+\infty)$,
$\gamma\in(0,+\infty)$ and~$C\in(0,+\infty)$ such that, if~$u$ is as in~\eqref{IVP},
then
\begin{equation}\label{SA}
\| u\|_{L^s(\Omega)}^{s-1+\gamma}(t)\le C\,\int_\Omega u^{s-1}(x,t)\,\NN[u](x,t)\,dx,
\end{equation}
where we used the notation
$$ \| u\|_{L^s(\Omega)}(t):=\left( {\int_\Omega u^s(x,t)\,dx}\right)^{1/s}.$$
For simplicity, we considered smooth solutions of~\eqref{IVP}
(in concrete cases, the notion of weak solutions may be treated similarly, see e.g.
the regularization methods discussed on page~235 of~\cite{VER}).
\medskip

After providing a general result on the decay of the solutions of~\eqref{IVP},
we will specify the operator~$\NN$ to the following concrete cases:
\begin{itemize}
\item the case of the Laplacian,
\item the case of the $p$-Laplacian,
\item the case of the porous medium equation,
\item the case of the doubly nonlinear equation,
\item the case of the mean curvature equation,
\item the case of the fractional Laplacian,
\item the case of the fractional $p$-Laplacian,
\item the sum of different space-fractional operators,
\item the case of the fractional porous medium equation,
\item the case of the fractional mean curvature equation.
\end{itemize}
The general result will be obtained by energy methods for nonlinear operators
(see e.g.~\cite{HAND}). Our approach will largely exploit a very deep and
detailed analysis of the time-fractional
evolution problems recently performed in~\cites{VER, VER2}
and, in a sense, our results can be seen as a generalization of those in~\cites{VER, VER2, ZZZ}
to comprise cases arising from space-fractional equations, nonlinear nonlocal operators,
geometric operators and nonlocal porous medium equations.\medskip

Also, the general framework that we provide can be useful to give a unified
setting in terms of energy inequalities.\medskip

Our ``abstract'' result is the following:

\begin{theorem}\label{THMDEC}
Let~$u$ be as in~\eqref{IVP}, under the structural condition in~\eqref{SA}.
Then,
\begin{equation} \label{v3}
\partial^\alpha_t \| u\|_{L^s(\Omega)}(t)\le-\frac{\| u\|_{L^s(\Omega)}^{\gamma}(t)}{C}.\end{equation}
Furthermore, 
\begin{equation} \label{v4}
\| u\|_{L^s(\Omega)}(t)\le\frac{C_\star}{1+t^{\alpha/\gamma}},
\end{equation}
for some~$C_\star>0$, possibly depending on~$C$, $\gamma$, $\alpha$ and~$\| u_0\|_{L^s(\Omega)}$.
\end{theorem}

We point out that the result in~\eqref{v4} is quite different from
the decay estimates for classical time-diffusion (compare e.g. with~\cite{HAND}).

Indeed, in~\eqref{v4}, a power-law decay is provided, while the classical uniformly elliptic
time-diffusion
case presents exponential decays. The power-law decay can be understood
by looking at the solution of
$$ \partial_t^\alpha e(t)=-e(t)$$
for~$t\in(0,+\infty)$ with initial datum~$e(0)=1$
and at the first Dirichlet eigenfunction~$\phi$ of a ball~$B$ normalized in such a way that
the corresponding
eigenvalue is equal to~$1$, namely
$$\left\{\begin{matrix} \Delta\phi=-\phi &{\mbox{ in }}B,\\
\phi=0&{\mbox{ on }}\partial B,\\
\|\phi\|_{L^2(B)}=1.&
\end{matrix}\right.$$
Then, the function~$u(x,t):=e(t)\,\phi(x)$ satisfies the fractional heat equation~$\partial^\alpha_t u=\Delta u$
in~$B$, with zero Dirichlet datum, and
\begin{equation} \label{909}\|u\|_{L^2(B)}= |e(t)|.\end{equation}
The function~$e$ is explicit in terms of the
Mittag-Leffler function (see e.g.~\cite{MAI, PAR} and the references therein),
and it satisfies~$e(t)\sim\frac{1}{t^\alpha}$ as~$t\to+\infty$.

This fact and~\eqref{909} imply a polynomial decay of the~$L^2$-norm of the solution,
in agreement with~\eqref{v4}.

The decay presented in~\eqref{v4} is also different from the case
of fast nonlinear diffusion, in which the solution
gets extincted in finite time, see e.g. Theorem~17 in~\cite{HAND}.
\medskip

We now specify Theorem~\ref{THMDEC} to several concrete cases,
also recovering the main results in~\cites{VER, ZZZ} and providing new applications.
Several new applications will be also given in~\cite{ELISA}.

\subsection{The cases of the Laplacian, of the $p$-Laplacian,
of the porous medium equation and of the doubly nonlinear equation}\label{RECOVER}

The doubly nonlinear operator (see e.g.~\cite{RAV}) is a general operator of the form
\begin{equation}\label{DOUB} u\longmapsto \Delta_p u^m,\end{equation}
with~$m\in(0,+\infty)$ and~$p\in(1,+\infty)$.

When~$m=1$, this operator reduces to the $p$-Laplacian
$$ \Delta_pu:={\rm div}\big( |\nabla u|^{p-2}\nabla u\big),$$
which in turn reduces to the classical Laplacian as~$p=2$.

When~$p=2$, the operator in~\eqref{DOUB} reduces to the porous medium operator (see e.g.~\cite{VAZ}
and the references therein)
\[ u\longmapsto \Delta u^m,\]
which again reduces to the Laplacian when~$m=1$.

In this setting, we have the following result:

\begin{theorem}\label{DBT}
Suppose that~$u$ is a solution of
$$ \partial^\alpha_t u = \Delta_p u^m$$
in~$\Omega\times(0,+\infty)$, with~$u(x,t)=0$ for
any~$x\in\partial\Omega$ and any~$t\geq0$.
Then, for any~$s\in(1,+\infty)$,
\begin{equation}\label{82394:2}
\| u\|_{L^s(\Omega)}(t)\le\frac{C}{1+t^{\frac\alpha{m(p-1)}}},
\end{equation}
for some~$C>0$.
\end{theorem}

As special cases of Theorem~\ref{DBT}, we can take~$m:=1$ and~$p:=2$,
which correspond to the $p$-Laplacian case and to the porous medium case, respectively.
We state these results explicitly for the convenience of the reader.

\begin{corollary}\label{DBT:1}
Suppose that~$u$ is a solution of
$$ \partial^\alpha_t u = \Delta_p u$$
in~$\Omega\times(0,+\infty)$, with~$u(x,t)=0$ 
for any~$x\in\partial\Omega$
and any~$t\geq0$.
Then, for any~$s\in(1,+\infty)$,
\[
\| u\|_{L^s(\Omega)}(t)\le\frac{C}{1+t^{\alpha/{(p-1)}}},
\]
for some~$C>0$.
\end{corollary}

\begin{corollary}\label{DBT:2}
Suppose that~$u$ is a solution of
$$ \partial^\alpha_t u = \Delta u^m$$
in~$\Omega\times(0,+\infty)$, with~$u(x,t)=0$ for
any~$x\in\partial\Omega$
and any~$t\geq0$.
Then, for any~$s\in(1,+\infty)$,
\[
\| u\|_{L^s(\Omega)}(t)\le\frac{C}{1+t^{\alpha/m}},
\]
for some~$C>0$.
\end{corollary}

When~$p\in(2,+\infty)$ in Corollary~\ref{DBT:1}, the case~$\alpha\nearrow1$ recovers the classical decay,
see e.g. Theorem 21 in~\cite{HAND}.\medskip

For results related to Corollary~\ref{DBT:1} when~$p=2$, see~\cites{RE1, RE2}.
Corollaries~\ref{DBT:1} and~\ref{DBT:2} can be compared with Theorems~8.1 and~9.1
in~\cite{VER}, respectively.

\subsection{The case of the mean curvature equation}

The setting in Theorem~\ref{THMDEC} is general enough
to deal with nonlinear operators of mean curvature type,
and to consider equations of the type
\begin{equation}\label{ZX1XGBOU}
\partial^\alpha_t u(x,t)={\rm div}\left( \frac{\nabla u(x,t)}{\sqrt{1+|\nabla u(x,t)|^2}}\right)
.\end{equation}
We recall that the right-hand-side of~\eqref{ZX1XGBOU} corresponds to
the mean curvature of the hypersurface described by the graph of the function~$u$,
see e.g. formula~(13.1) in~\cite{GIU}. The result obtained in this setting goes as follows:

\begin{theorem}\label{MCE}
Suppose that~$u$ is a solution of~\eqref{ZX1XGBOU}
in~$\Omega\times(0,+\infty)$, with~$u(x,t)=0$
for any~$x\in\partial\Omega$
and any~$t\geq0$.
Assume that either
\begin{equation}\label{GBOU}
n\in\{1,2\}\quad{\mbox{ and }}\quad
\sup_{t>0} \int_\Omega
\sqrt{1+|\nabla u(x,t)|^2}\,dx<+\infty
\end{equation}
or
\begin{equation}\label{GBOU2}
\sup_{x\in\Omega,\;t>0} |\nabla u(x,t)| <+\infty.
\end{equation}
Then, for any~$s\in(1,+\infty)$,
\begin{equation}\label{SPE} \| u\|_{L^s(\Omega)}(t)\le\frac{C}{1+t^{\alpha}},\end{equation}
for some~$C>0$.\end{theorem}

\subsection{The case of the fractional Laplacian,
of the fractional $p$-Laplacian, and of the sum of different space-fractional operators}

The setting in Theorem~\ref{THMDEC} is general enough to comprise also the case
of operators modeling spatial nonlocal diffusion of fractional kind. The main example of
such operators is given by the fractional Laplacian of order~$\sigma\in(0,1)$,
which can be defined (up to a multiplicative constant that we neglect for simplicity) by
\begin{equation}\label{sLA}
(-\Delta)^\sigma u(x):=\int_{\R^n} \frac{u(x)-u(y)}{|x-y|^{n+2\sigma}}\,dy.
\end{equation}
Here and in the following, we 
implicitly suppose that these types of singular integrals are taken in the principal value sense.
The fractional Laplacian provides a natural framework for many problems in theoretical and
applied mathematics, see e.g.~\cites{LAN, SIL, BUC} and the references therein.

Several nonlinear variations of the fractional Laplacian can be taken into account,
see e.g.~\cite{KUU0, KUU, PUCCI, 7UOA11SQ, COZZ}
and the references therein. In particular, for any~$p\in(1,+\infty)$, one can consider the operator
\begin{equation}\label{psLA}
(-\Delta)_p^\sigma u(x):=\int_{\R^n} \frac{\big|
u(x)-u(y)\big|^{p-2}\big(u(x)-u(y)\big)}{|x-y|^{n+\sigma p}}\,dy.
\end{equation}
Of course, when~$p=2$ the operator in~\eqref{psLA} reduces to that in~\eqref{sLA}.
In this setting, we have the following decay result for solutions of fractional-time equations
whose spatial diffusion is driven by the nonlinear fractional
operator in~\eqref{psLA}.

\begin{theorem}\label{THp1}
Suppose that~$u$ is a solution of
$$ \partial^\alpha_t u(x,t)+
(-\Delta)^\sigma_p u(x,t)=0$$
in~$\Omega\times(0,+\infty)$, with~$u(x,t)=0$ for
any~$x\in\R^n\setminus\Omega$ and
for any~$t>0$.
Then, for any~$s\in(1,+\infty)$,
\begin{equation}\label{PERO2}
\| u\|_{L^s(\Omega)}(t)\le\frac{C}{1+t^{\alpha/(p-1)}},
\end{equation}
for some~$C>0$.
\end{theorem}

When~$p=2$, decay estimates for nonlocal equations
have been very recently obtained in~\cite{ZZZ}.\medskip

An extension of Theorem~\ref{THp1} holds true also for sums of different,
possibly nonlinear, space-fractional diffusion operators:

\begin{theorem}\label{THp2}
Let~$N\in\N$, $N\ge1$. Let~$\sigma_1,\dots,\sigma_N\in(0,1)$,
$p_1,\dots,p_N\in(1,+\infty)$
and~$\beta_1,\dots,\beta_N\in(0,+\infty)$.
Suppose that~$u$ is a solution of
$$ \partial^\alpha_t u(x,t)+\sum_{j=1}^N
\beta_j\,(-\Delta)^{\sigma_j}_{p_j} u(x,t)=0$$
in~$\Omega\times(0,+\infty)$, with~$u(x,t)=0$ for any~$x\in\R^n\setminus\Omega$ and
for any~$t>0$.
Then, for any~$s\in(1,+\infty)$,
\begin{equation}\label{PERO22}
\| u\|_{L^s(\Omega)}(t)\le\frac{C}{1+t^{\alpha/(p_{\max}-1)}},\end{equation}
for some~$C>0$,
where
$$ p_{\max}:=\max\{p_1,\dots,p_N\}.$$
\end{theorem}

More general settings for superpositions of fractional
operators can be also considered in our setting, see e.g.~\cite{CABSE}.
\medskip

Another interesting case arises
from the sum of
fractional operators in different directions.
Precisely, fixed~$j\in\{1,\dots,n\}$ one can consider the unit vector~$e_j$
(i.e., the $j$th element of
the Euclidean basis of~$\R^n$), and define
the fractional
Laplacian in direction~$e_j$, namely
$$ (-\partial^2_{x_j})^{\sigma_j} u(x):=
\int_\R \frac{u(x)-u(x+\rho e_j)}{\rho^{1+2\sigma_j}}\,d\rho,$$
with~$\sigma_j\in(0,1)$.
Then, given~$\beta_1,\dots,\beta_n>0$, one can consider
the superposition of such operators, that is
$$ (-\Delta_{\beta})^\sigma u(x):=\sum_{j=1}^n \beta_j
(-\partial^2_{x_j})^{\sigma_j} u(x).$$
Here, we are using the formal notation~$\sigma:=(\sigma_1,\dots,\sigma_n)$ and~$\beta:=(\beta_1,\dots,\beta_n)$.
Notice that~$(-\Delta_{\beta})^\sigma$ is similar
to,
but structurally very different
from, the fractional Laplacian,
since the nonlocal character of the fractional Laplacian
also takes into account the
interactions in directions different than~$e_1,\dots,e_n$:
for instance, even if~$
\sigma_1=\dots=\sigma_n=1/2$ and~$\beta_1=\dots=\beta_n=1$,
the operator~$(-\Delta_{\beta})^\sigma$ does not
reduce to the square root of the Laplacian.

These types of anisotropic fractional operators (and even more general ones)
have been considered in~\cites{ROS1, ROS2, ROS3, FVX}.
In our setting, we have the following decay result:

\begin{theorem}\label{THp3}
Suppose that~$u$ is a solution of
$$ \partial^\alpha_t u(x,t)+(-\Delta_{\beta})^\sigma u(x,t)=0$$
in~$\Omega\times(0,+\infty)$, with~$u(x,t)=0$ for
any~$x\in\R^n\setminus\Omega$ and
for any~$t>0$.
Then, for any~$s\in(1,+\infty)$,
\begin{equation*}
\| u\|_{L^s(\Omega)}(t)\le\frac{C}{1+t^{\alpha}},\end{equation*}
for some~$C>0$.
\end{theorem}

More general operators, such as the sum of fractional Laplacians
along linear subspaces of~$\R^n$, as well as
operators in integral superposition, may also be taken into
account in Theorem~\ref{THp3}, but we focused on an explicit
case for simplicity of notations.

\subsection{The case of the fractional porous medium equation}

We consider here a porous medium diffusion operator of fractional
type, given by
$$ u\longmapsto (-\Delta)^\sigma u^m,$$
with~$\sigma\in(0,1)$ and~$m\in(0,+\infty)$, where~$(-\Delta)^\sigma$
is the fractional Laplace operator defined in~\eqref{sLA}.

In the classical time-diffusion case, such equation has been introduced
and analyzed in~\cite{PM1, PM2} (remarkably, in this case,
any nontrivial nonnegative solution becomes strictly positive instantaneously,
and this is a different feature with respect to the classical
porous medium equation).\medskip

In our setting, we will consider the time-fractional
version of the space-fractional porous medium equation
and establish the following decay estimate:

\begin{theorem}\label{VVV}
Suppose that~$u$ is a solution of
$$ \partial^\alpha_t u +(-\Delta)^\sigma u^m=0$$
in~$\Omega\times(0,+\infty)$, 
with~$u(x,t)=0$ for any~$x\in
\R^n\setminus\Omega$ and any~$t\ge0$.
Then, for any~$s\in(1,+\infty)$,
\begin{equation}\label{sdfbgssh82394:2}
\| u\|_{L^s(\Omega)}(t)\le\frac{C}{1+t^{\frac\alpha{m}}},
\end{equation}
for some~$C>0$.
\end{theorem}

Similar results may also be obtained in more general settings
for doubly nonlinear and doubly fractional porous medium
equations. We also observe that, for~$m=1$,
Theorem~\ref{VVV} boils down to Theorem~\ref{THp1}
with~$p:=2$. 

\subsection{The case of the fractional mean curvature equation}

The notion of nonlocal perimeter functional has been introduced
and analyzed in~\cite{CRS}. While the first variation of the classical
perimeter functional consists in the mean curvature operator,
the first variation of the nonlocal perimeter produces
an object, which can be seen as a nonlocal mean curvature,
and which corresponds to a weighted average of the characteristic
function of a set with respect to a singular kernel.
The study of such fractional mean curvature operator
is a very interesting topic of
research in itself and the recent literature produced
several contributions in this context,
see e.g.~\cites{CCN1, CCN2, CCN2a, CCN2b, CCN3, CCN3a, CCN3b, CCN4}.
The nonlocal mean curvature also induces a geometric flow,
as studied in~\cite{PONS, MARIEL, MATTEO, SINES}.
See also~\cite{123} for a recent survey on the topic of nonlocal
minimal surfaces and nonlocal mean curvature equations.\medskip

For smooth hypersurfaces with a structure of complete graphs, the notion of nonlocal mean curvature
can be introduced as follows (see e.g. formula~(3.5) in~\cite{CCN3a}).
For any~$r\in\R$ and~$\sigma\in(0,1)$, we set
\begin{equation}\label{F} F(r):=\int_0^r \frac{d\tau}{(1+\tau^2)^{(n+1+\sigma)/2}}\end{equation}
and we consider the (minus) nonlocal mean curvature operator corresponding to the choice
\begin{equation}\label{CH6}
{\mathcal{N}}[u](x,t):=\int_{\R^n} \frac{1}{|y|^{n+\sigma}}
\,F\left( \frac{u(x,t)-u(x+y,t)}{|y|}\right)\,dy.\end{equation}
In this setting, we provide a decay estimate
for graphical solutions of the fractional mean curvature equation,
as stated in the following result:

\begin{theorem}\label{FF}
Suppose that~$u$ is a solution of
$$ \partial^\alpha_t u(x,t)=
\int_{\R^n} \frac{1}{|y|^{n+\sigma}}
\,F\left( \frac{u(x+y,t)-u(x,t)}{|y|}\right)\,dy$$
with~$u(x,t)=0$ for any~$x\in\R^n\setminus\Omega$ and
for any~$t>0$.

Assume that
\begin{equation}\label{GBOU3}
\sup_{x\in\Omega,\;t>0} |\nabla u(x,t)| <+\infty.
\end{equation}
Then, for any~$s\in(1,+\infty)$,
\begin{equation*} \| u\|_{L^s(\Omega)}(t)\le\frac{C}{1+t^{\alpha}},\end{equation*}
for some~$C>0$.
\end{theorem}

The recent literature has considered the evolution
of graphs under the fractional mean curvature flow, see Section~6
of~\cite{MARIEL}. In this respect,
Theorem~\ref{FF} here can be seen as the first study of
evolution equations driven by the fractional mean curvature in which the flow
possesses a memory effect.\medskip

The rest of the paper is devoted to the proofs of the above mentioned results. First
we prove the general statement of
Theorem~\ref{THMDEC}, and then we check that condition~\eqref{SA}
is verified in all the concrete cases taken into consideration.

\section{Proofs}

We will exploit the following result,
which follows from Corollary~3.1 in~\cite{VER}.

\begin{lemma}\label{AZ}
Let~$s>1$, $u:\Omega\times[0,+\infty)\to\R$ and~$u_0(x):=u(x,0)$.
Let~$v(t):=\|u\|_{L^s(\Omega)}(t)$ and suppose that~$u_0\in L^s(\Omega)$
and, for every~$T>0$, that~$u\in L^s((0,T),L^s(\Omega))$. Then,
$$ v^{s-1}(t)\,\partial^\alpha_t v(t)\le \int_\Omega u^{s-1}(x,t)\,\partial^\alpha_t u(x,t)\,dx.$$
\end{lemma}

\subsection{Proof of Theorem~\ref{THMDEC}}

Without loss of generality, we can suppose that~$\|u_0\|_{L^s(\Omega)}\ne0$,
and we set~$v(t):=\|u\|_{L^s(\Omega)}(t)$.
Hence, recalling~\eqref{IVP} and Lemma~\ref{AZ},
\begin{equation}\label{21a}
v^{s-1}(t)\,\partial^\alpha_t v(t)\le -\int_\Omega u^{s-1}(x,t)\,\NN[ u](x,t)\,dx.\end{equation}
Using this and~\eqref{SA}, we thus find that
\begin{equation} \label{v1}
v^{s-1}(t)\,\partial^\alpha_t v(t)\le-\frac{1}{C}\,\| u\|_{L^s(\Omega)}^{s-1+\gamma}(t)
=-\frac{v^{s-1+\gamma}(t)}{C}.\end{equation}
{F}rom~\eqref{v1}, we plainly obtain~\eqref{v3}
at all~$t$ for which~$v(t)\ne0$.

But at the points~$t$ at which~$v(t)=0$,
we see that~\eqref{v3} is also automatically\footnote{For another
approach allowing for the division by the prefactor in~\eqref{v1},
see Lemma~2.1 in the recent preprint~\cite{DP}.}
satisfied in view of
the following observation: using that~$u$ is smooth,
combined with the H\"older's Inequality with exponents~$s/(s-1)$ and~$s$,
we have that, a.e.~$t>0$,
\begin{eqnarray*}
\left|\frac{\partial}{\partial t} v(x,t)\right|&=&
\left|\frac{\partial}{\partial t} \left( \int_\Omega |u(x,t)|^s\,dx\right)^{1/s}\right|
\\&\le&\left( \int_\Omega |u(x,t)|^s\,dx\right)^{(1-s)/s}
\,\int_\Omega |u(x,t)|^{s-1} \,\left| \frac{\partial u}{\partial t}(x,t)\right|\,dx
\\&\le& \left( \int_\Omega |u(x,t)|^s\,dx\right)^{(1-s)/s}
\left( \int_\Omega |u(x,t)|^s\,dx\right)^{(s-1)/s}\left(\int_\Omega
\left| \frac{\partial u}{\partial t}(x,t)\right|^s\,dx\right)^{1/s}\\
&=& \left\| \frac{\partial u}{\partial t}\right\|_{L^s(\Omega)}(t).
\end{eqnarray*}
Hence, $v$ is
Lipschitz continuous, and therefore
we see that
$$ \lim_{\tau\nearrow t} \frac{v(t)-v(\tau)}{(t-\tau)^\alpha}=0,$$
and, as a consequence,
\begin{eqnarray*}
\frac{d}{dt}\int_0^t \frac{v(t)-v(\tau)}{(t-\tau)^\alpha}\,d\tau&=&
\int_0^t \frac{\partial_t v(t)}{(t-\tau)^\alpha}\,d\tau
-\alpha\int_0^t \frac{v(t)-v(\tau)}{(t-\tau)^{1+\alpha}}\,d\tau\\
&=&
\frac{\partial_t v(t) \,t^{1-\alpha} }{1-\alpha}
-\alpha\int_0^t \frac{v(t)-v(\tau)}{(t-\tau)^{1+\alpha}}\,d\tau
.\end{eqnarray*}
Comparing this with~\eqref{K-DER}, we find that
\begin{eqnarray*}
\partial_t^{\alpha} v(t)&=&
\frac{d}{dt}
\int_0^t \frac{ v(\tau)-v(t)}{ (t-\tau)^\alpha }\,d\tau
+
\frac{d}{dt}
\int_0^t \frac{ v(t)-v(0)}{ (t-\tau)^\alpha }\,d\tau
\\&=& -\frac{\partial_t v(t) \,t^{1-\alpha} }{1-\alpha}
+\alpha\int_0^t \frac{v(t)-v(\tau)}{(t-\tau)^{1+\alpha}}\,d\tau+
\frac{d}{dt}\frac{\big(v(t)-v(0)\big)\,t^{1-\alpha}}{1-\alpha}\\&=&
-\frac{\partial_t v(t) \,t^{1-\alpha} }{1-\alpha}
+\alpha\int_0^t \frac{v(t)-v(\tau)}{(t-\tau)^{1+\alpha}}\,d\tau+
\frac{\partial_tv(t)\,t^{1-\alpha}}{1-\alpha}
+\frac{v(t)-v(0)}{t^\alpha}\\
&=&\alpha\int_0^t \frac{v(t)-v(\tau)}{(t-\tau)^{1+\alpha}}\,d\tau
+\frac{v(t)-v(0)}{t^\alpha}.
\end{eqnarray*}
Therefore, at points~$t$ where~$v(t)=0$, using that~$v\ge0$,
we see that
$$ \partial_t^{\alpha} v(t)=-\alpha\int_0^t \frac{
v(\tau)}{(t-\tau)^{1+\alpha}}\,d\tau
-\frac{v(0)}{t^\alpha}\le0,$$
which gives that~\eqref{v3} is satisfied in this case as well.

Now, we
prove~\eqref{v4}. To this aim, we consider
the solution~$w(t)$ of the nonlinear fractional differential equation 
\begin{equation}\label{w2} \left\{\begin{matrix}
\partial^\alpha_t w(t) =-\displaystyle\frac{w^{\gamma}(t)}{C} & {\mbox{ for any }}t>0,\\ \\
w(0)=v(0).
\end{matrix}\right.\end{equation}
When~$\gamma=1$, the
function~$w$ is explicitly known in terms of
the Mittag-Leffler function, see~\cites{MAI, PAR}. The general case~$\gamma>0$
has been dealt with in detail in Section~7 of~\cite{VER}. In particular (see Theorem~7.1
in~\cite{VER}) it holds that
\begin{equation}\label{w}
w(t)\le \frac{C_\star}{1+t^{\alpha/\gamma}},
\end{equation}
for some~$C_\star>0$, possibly depending on~$C$, $\gamma$, $\alpha$ and~$v(0)$.
Moreover, by \eqref{2.2bis}, \eqref{w2} and the
comparison principle (see e.g. Lemma~2.6 in~\cite{VER}), 
we have that~$v(t)\leq w(t)$ for all~$t\ge0$.
Using this
and~\eqref{w}, we obtain~\eqref{v4}.~\hfill$\Box$

\begin{remark}{\rm We observe that the constant~$C$ in~\eqref{v3}
is exactly the one coming from~\eqref{SA}. If needed, the long time
behavior in~\eqref{v4} can be also made more
precise in terms of $\| u_0\|_{L^s(\Omega)}$.
Indeed, recalling formula~(41) in~\cite{VER},
in the notation used for the proof of Theorem~\ref{THMDEC}
one can define
$$ t_0:=\bar{C}\,w^{(1-\gamma)/\alpha}(0),$$
with~$\bar{C}>0$ depending
only on~$\alpha$, $\gamma$ and~$C$, and
$$ \bar w(t):=\begin{cases}
w(0) & {\mbox{ if }}t\in[0,t_0],\\
w(0)\,(t_0/t)^{\alpha/\gamma}& {\mbox{ if }}t\in(t_0,+\infty),
\end{cases}$$
and conclude that~$w(t)\le \bar w(t)$. In this way, for large $t$, we have that
$$ \|u\|_{L^s(\Omega)}(t)=v(t)\le\bar w(t)=
\frac{w(0)\,t_0^{\alpha/\gamma}}{t^{\alpha/\gamma}}=
\frac{\|u_0\|_{L^s(\Omega)}\,
\big( \bar{C}\,
\|u_0\|_{L^s(\Omega)}^{(1-\gamma)/\alpha}
\big)^{\alpha/\gamma}}{t^{\alpha/\gamma}}
=
\frac{
{\tilde{C}}\,
\|u_0\|_{L^s(\Omega)}^{1/\gamma}}{t^{\alpha/\gamma}}
,$$
with~$\tilde{C}>0$
depending
only on~$\alpha$, $\gamma$ and~$C$.
}\end{remark}

\subsection{Proof of Theorem~\ref{DBT}}

We set
\begin{equation}\label{09sbala}
v:=u^{\frac{s-2+p+(m-1)(p-1)}{p}}\end{equation} and we point out that
\begin{equation}\label{72etr9ig23fgerkhgg}
\begin{split} |\nabla v|^p \,&=\left( \frac{s-2+p+(m-1)(p-1)}{p}\right)^p
u^{{s-2+(m-1)(p-1)}}\,|\nabla u|^p\\&=
\left( \frac{s-2+p+(m-1)(p-1)}{p}\right)^p
u^{s-2}\nabla u \cdot\big( u^{(m-1)(p-1)}|\nabla u|^{p-2}\nabla u\big)\\&
=\left( \frac{s-2+p+(m-1)(p-1)}{p}\right)^p
\,\frac{1}{(s-1)\,m^{p-1}}
\,\nabla u^{s-1} \cdot\big( |\nabla u^m|^{p-2}\nabla u^m\big).\end{split}\end{equation}
Now, when~$p\in(1,n)$, we recall the Sobolev exponent
$$ p_\star:=\frac{np}{n-p}$$
and we claim that, if~$p\in(1,n)$ and~$q\in[1,p_\star]$,
as well as if~$p\in[n,+\infty)$ and~$q\in[1,+\infty)$, it holds that 
\begin{equation}\label{CLAGA}
\| v\|_{L^q(\Omega)}^p(t)\le C_0\,\int_{\Omega} |\nabla v(x,t)|^p\,dx,
\end{equation}
for some~$C_0>0$. Indeed, when~$p\in(1,n]$, the inequality
in~\eqref{CLAGA} follows from the Sobolev Embedding Theorem.
When instead~$p\in(n,+\infty)$ we can use the
Sobolev Embedding Theorem with exponent~$n$
to write
$$ \| v\|_{L^q(\Omega)}^n(t)\le C_0\,\int_{\Omega} |\nabla v(x,t)|^n\,dx.$$
Combining this
with the H\"older's Inequality for the norm of the gradient, we obtain~\eqref{CLAGA}
(up to renaming constants).

We also observe that when~$p\in(1,n)$ 
and
\begin{equation}\label{COA}
s\ge \max\left\{m-\frac{1}{p-1},\;
\frac{n\big(1-m(p-1) \big)}p
\right\},
\end{equation}
it holds that
\begin{equation}\label{RANGO}
\frac{sp}{s-2+p+(m-1)(p-1)}\in [1,p_\star].
\end{equation}
Indeed, we have that
$$s-2+p+(m-1)(p-1)>1-2+p+(m-1)(p-1)=m(p-1)\ge0.$$
Moreover,
$$ s-2+p+(m-1)(p-1)-sp = m(p-1)-1-s(p-1)\le 0$$
thanks to~\eqref{COA}. This gives that~$\frac{sp}{s-2+p+(m-1)(p-1)}\ge1$.

In addition,
\begin{eqnarray*}
s(n-p)-n\big( s-2+p+(m-1)(p-1) \big) =
-sp-n\big(p-2+(m-1)(p-1) \big)\le0,
\end{eqnarray*}
due to~\eqref{COA}, which gives that~$\frac{sp}{s-2+p+(m-1)(p-1)}\le p_\star$.
These considerations prove~\eqref{RANGO}.

By means of~\eqref{RANGO}, when either~$p\in[n,+\infty)$ or~\eqref{COA} holds true,
we can choose~$q:=\frac{sp}{s-2+p+(m-1)(p-1)}$ in~\eqref{CLAGA}.
Hence, recalling~\eqref{09sbala}, we find that
\begin{eqnarray*}
\|u\|_{L^s(\Omega)}^{
{s-2+p+(m-1)(p-1)} }(t)&=&
\left(\int_\Omega 
u^{s}(x,t)\,dx\right)^{
\frac{s-2+p+(m-1)(p-1)}{s}}\\&
=&
\left(\int_\Omega 
u^{\frac{s-2+p+(m-1)(p-1)}{p}\cdot\frac{sp}{s-2+p+(m-1)(p-1)}}(x,t)\,dx\right)^{
\frac{s-2+p+(m-1)(p-1)}{s}}\\&
=&
\left(\int_\Omega v^{\frac{sp}{s-2+p+(m-1)(p-1)}}(x,t)\,dx\right)^{
\frac{s-2+p+(m-1)(p-1)}{s}}
\\&=&
\| v\|_{L^{\frac{sp}{s-2+p+(m-1)(p-1)}}(\Omega)}^p(t)
\\&\le& C_0\,\int_{\Omega} |\nabla v(x,t)|^p\,dx.
\end{eqnarray*}
As a consequence, making use of~\eqref{72etr9ig23fgerkhgg},
we conclude that
$$ \|u\|_{L^s(\Omega)}^{
{s-2+p+(m-1)(p-1)} }(t)\le C_1\int_\Omega
\nabla u^{s-1} \cdot\big( |\nabla u^m|^{p-2}\nabla u^m\big)\,dx,$$
provided that either~$p\in[n,+\infty)$ or~\eqref{COA} holds true.

This gives that condition~\eqref{SA}
is satisfied in this case with~$\gamma:=m(p-1)$. This and~\eqref{v4}
imply that, if either~$p\in[n,+\infty)$ or~\eqref{COA} holds true, then
\begin{equation}\label{82394}
\| u\|_{L^s(\Omega)}(t)\le\frac{C_\star}{1+t^{\frac\alpha{m(p-1)}}},
\end{equation}
for some~$C_\star>0$.
Also, when~\eqref{COA} is not satisfied, we have that
$$s< \max\left\{m-\frac{1}{p-1},\;
\frac{n\big(1-m(p-1) \big)}p
\right\}=:\bar s$$
and in this case the H\"older's Inequality implies
that~$\| u\|_{L^s(\Omega)}(t)\le C\,\| u\|_{L^{\bar s}(\Omega)}(t)$,
for some~$C>0$,
and~$\bar s$ lies in the range satisfying~\eqref{82394}.

This observation and~\eqref{82394} imply~\eqref{82394:2}, as desired.~\hfill$\Box$

\subsection{Proof of Theorem~\ref{MCE}}

We set~$v:=u^{s/2}$. Notice that
\begin{equation}\label{SOC:1}
|\nabla v|^2 = \frac{s^2}{4}\,u^{s-2}|\nabla u|^2=
\frac{s^2}{4(s-1)}\,\nabla u\cdot\nabla u^{s-1}.
\end{equation}
We distinguish two cases, according to~\eqref{GBOU} and~\eqref{GBOU2}.
We first consider the case in which~\eqref{GBOU}
holds true. Then, by Cauchy-Schwarz Inequality,
\begin{equation}\label{SOC:2}\begin{split}
\int_\Omega |\nabla v(x,t)|\,dx\,&=
\int_\Omega \frac{|\nabla v(x,t)|}{\big(1+|\nabla u(x,t)|^2\big)^{1/4}}
\,\big(1+|\nabla u(x,t)|^2\big)^{1/4}
\,dx\\ &\le \sqrt{
\int_\Omega \frac{|\nabla v(x,t)|^2}{\sqrt{1+|\nabla u(x,t)|^2}}
\,dx}\,\sqrt{
\int_\Omega 
\sqrt{1+|\nabla u(x,t)|^2}
\,dx}.\end{split}
\end{equation}
Moreover, when~$n>1$, from the 
Gagliardo-Nirenberg-Sobolev Inequality, we know that
$$ \left(\int_\Omega u^{\frac{sn}{2(n-1)}}(x,t)\,dx\right)^{\frac{n-1}n}=
\| v\|_{ L^{\frac{n}{n-1} } (\R^n)} (t)
\le C_0\int_{\R^n} |\nabla v(x,t)|\,dx$$
for some~$C_0>0$.
Also, when~$n=1$, one can use the Fundamental Theorem of Calculus and check that,
for any~$q\in[1,+\infty)$,
$$ \left(\int_\Omega u^{\frac{sq}{2}}(x,t)\,dx\right)^{\frac{1}q}=
\| v\|_{ L^{q } (\R^n)} (t)
\le C_1\int_{\R^n} |\nabla v(x,t)|\,dx$$
for some~$C_1>0$.

Using this, \eqref{SOC:1} and~\eqref{SOC:2}, we obtain that
$$ \left(\int_\Omega u^{\frac{sq}{2}}(x,t)\,dx\right)^{\frac{1}q}\le
C_2
\sqrt{
\int_\Omega \frac{\nabla u(x,t)\cdot\nabla u^{s-1}(x,t)
}{\sqrt{1+|\nabla u(x,t)|^2}}
\,dx}\,\sqrt{
\int_\Omega
\sqrt{1+|\nabla u(x,t)|^2}
\,dx},$$
where~$q=\frac{n}{n-1}$ when~$n=2$, and any~$q\in[1,+\infty)$ when~$n=1$.
{F}rom this and assumption~\eqref{GBOU}, we find that
$$ \left(\int_\Omega u^{\frac{sq}{2}}(x,t)\,dx\right)^{\frac{2}q}\le
C_3
\,\int_\Omega \frac{\nabla u(x,t)\cdot\nabla u^{s-1}(x,t)
}{\sqrt{1+|\nabla u(x,t)|^2}}
\,dx,$$
where~$q=\frac{n}{n-1}$ when~$n=2$, and any~$q\in[1,+\infty)$ when~$n=1$.
In any case, when~$n\in\{1,2\}$, we have that
we can take~$q=2$ and write
$$ \int_\Omega u^{s}(x,t)\,dx\le
C_3
\,\int_\Omega \frac{\nabla u(x,t)\cdot\nabla u^{s-1}(x,t)
}{\sqrt{1+|\nabla u(x,t)|^2}}
\,dx.$$
Therefore, we have that~\eqref{SA}
is satisfied for any~$s\in(1,+\infty)$ and~$\gamma:=1$.
This and \eqref{v4} imply \eqref{SPE}, as desired.

Now we deal with the case in which~\eqref{GBOU2} is satisfied.
We can also assume that~$n\ge3$ (since the cases~$n\in\{1,2\}$ are covered by~\eqref{GBOU}).
Then, exploiting the Gagliardo-Nirenberg-Sobolev Inequality in this situation
and recalling~\eqref{SOC:1}, we see that
\begin{eqnarray*}
\frac{s^2}{4(s-1)}\,
\int_\Omega \nabla u(x,t)\cdot\nabla u^{s-1}(x,t)\,dx&=&\int_\Omega
|\nabla v(x,t)|^2\,dx\\&\ge&C_0\,
\| v\|^2_{L^{\frac{2n}{n-2}}(\Omega)}\\
&=&C_0\,\left( \int_\Omega u^{{\frac{sn}{n-2}}}\right)^{\frac{n-2}{n}}
,
\end{eqnarray*}
for some~$C_0>0$. Hence, by H\"older's Inequality,
$$ \int_\Omega \nabla u(x,t)\cdot\nabla u^{s-1}(x,t)\,dx\ge C_1\,\|u\|_{L^s(\Omega)}^s,$$
for some~$C_1>0$. Thus, in light of~\eqref{GBOU2},
$$ \int_\Omega \frac{\nabla u(x,t)\cdot\nabla u^{s-1}(x,t)
}{\sqrt{1+|\nabla u(x,t)|^2}}
\,dx\ge C_2\,\|u\|_{L^s(\Omega)}^s.$$
This gives that~\eqref{SA} holds true in this case with~$\gamma:=1$.
Therefore, by means of~\eqref{v4} we obtain~\eqref{SPE}, as desired.~\hfill$\Box$

\subsection{Proof of Theorem~\ref{THp1}}

We define~$v:=u^{(s-2+p)/p}$ and we claim that
\begin{equation}\label{DISV}
|v(x,t)-v(y,t)|^p\le C_0 
|u(x,t)-u(y,t)|^{p-2}(u(x,t)-u(y,t))(u^{s-1}(x,t)-u^{s-1}(y,t)),
\end{equation}
for some~$C_0>0$. To this aim, we consider the auxiliary function
\begin{equation*} (1,+\infty)\ni\lambda\longmapsto g(\lambda):=
\frac{(\lambda^{(s-2+p)/p}-1)^p}{(\lambda-1)^{p-1} (\lambda^{s-1}-1)}.\end{equation*}
We recall that~$s-2+p>-1+p>0$ and observe that
$$ \lim_{\lambda\to+\infty} g(\lambda)=
\lim_{\lambda\to+\infty}
\frac{\left(1-\frac1{\lambda^{(s-2+p)/p}}\right)^p}{
\left(1-\frac1\lambda\right)^{p-1} \left(1-\frac1{\lambda^{s-1}}\right)}=
\frac{(1-0)^p}{(1-0)^{p-1} (1-0)}=1$$
and that
\begin{equation*}\begin{split}& \lim_{\lambda\searrow1} g(\lambda)=\lim_{\e\searrow0}
\frac{\big((1+\e)^{(s-2+p)/p}-1\big)^p}{\big((1+\e)-1\big)^{p-1} 
\big( (1+\e)^{s-1}-1\big)}\\&\qquad=
\lim_{\e\searrow0}
\frac{\big(1 + ((s-2+p)/p)\e
+o(\e)-1\big)^p}{\e^{p-1} 
\big( 1+(s-1)\e +o(\epsilon)-1\big)}=
\frac{((s-2+p)/p)^p }{s-1}.\end{split}\end{equation*}
Consequently,
$$ C_0:=\sup_{\lambda\in(1,+\infty)}g(\lambda) <+\infty.$$
Now, to prove~\eqref{DISV}, we may suppose,
up to exchanging~$x$ and~$y$,
that~$u(x,t)\geq u(y,t)$. 
Also, when either~$u(y,t)=0$
or~$u(x,t)=u(y,t)$, then~\eqref{DISV} is obvious.
Therefore, we can assume that~$u(x,t)>u(y,t)>0$
and set
$$ \lambda(x,y,t):=\frac{u(x,t)}{u(y,t)}\in(1,+\infty)$$
and conclude that
\begin{equation}\label{Com}
\begin{split}
C_0\, &\ge g(\lambda(x,y,t))\\
&=
\frac{\left(\frac{u^{(s-2+p)/p}(x,t)}{u^{(s-2+p)/p}(y,t)} -1\right)^p}{\left(
\frac{u(x,t)}{u(y,t)}-1\right)^{p-1}
\left(\frac{u^{s-1}(x,t)}{u^{s-1}(y,t)}-1\right)}
\\ &=
\frac{ \left( u^{(s-2+p)/p}(x,t)-u^{(s-2+p)/p}(y,t)\right)^p}{\left(
{u(x,t)}-{u(y,t)}\right)^{p-1}
\left( {u^{s-1}(x,t)}-{u^{s-1}(y,t)}\right)},\end{split}
\end{equation}
and this proves~\eqref{DISV}.

Now, when~$p\in\left(1,\frac{n}\sigma\right)$,
we consider the fractional critical exponent
\begin{equation}\label{psigma}p_\sigma:=
\frac{np}{n-\sigma p}.\end{equation}
We claim that
\begin{equation}\label{EM} 
\| v\|_{L^{q}(\Omega)}^{p}(t) \le C_1\,
\iint_{\R^{2n}} \frac{|v(x,t)-v(y,t)|^p}{
|x-y|^{n+\sigma p}}\,dx\,dy,\end{equation}
for some~$C_1>0$, for every~$q\in[1,p_\sigma]$
when~$p\in\left(1,\frac{n}\sigma\right)$, and for every~$q\in[1,+\infty)$ when~$p
\in\left[\frac{n}\sigma,+\infty\right)$.
Indeed,
when~$p\in\left(1,\frac{n}\sigma\right)$, then~\eqref{EM}
follows by the Gagliardo-Sobolev Embedding
(see e.g. Theorem~3.2.1 
in~\cite{BUC}). If instead~$
p\in 
\left[\frac{n}\sigma,+\infty\right)$ and~$q\in[1,+\infty)$,
we set~$\tilde q:=1+\max\{p,q\}$. Notice that
$$ 0<\frac{n}{p}-\frac{n}{\tilde q}<\frac{n}{p}\le\sigma<1.$$
Hence, we can take
$$ \tilde\sigma\in \left( \frac{n}{p}-\frac{n}{\tilde q},\;\frac{n}{p}\right)$$
and, since
$$p\in\left(1,\frac{n}{\tilde\sigma}\right)\quad{\mbox{
and }}\quad\tilde q\le \frac{np}{n-\tilde\sigma p}=p_{\tilde\sigma},$$
we can make use
of Gagliardo-Sobolev Embedding
(see e.g. Theorem~3.2.1 
in~\cite{BUC}) with exponents~$\tilde\sigma$
and~$\tilde q$. In this way, we find that
\begin{equation}\label{MAK0}
\| v\|_{L^{\tilde q}(\Omega)}^{p}(t) \le C_\sharp\,
\iint_{\R^{2n}} \frac{|v(x,t)-v(y,t)|^p}{
|x-y|^{n+\tilde\sigma p}}\,dx\,dy,\end{equation}
for some~$C_\sharp>0$. Now, we fix~$M>0$, to be
taken appropriately large, and we observe that
\begin{eqnarray*}
&& \iint_{\R^{2n}\cap\{ |x-y|\le M\}} \frac{|v(x,t)-v(y,t)|^p}{
|x-y|^{n+\tilde\sigma p}}\,dx\,dy\le
M^{(\sigma-\tilde\sigma)p}\iint_{\R^{2n}\cap\{ |x-y|\le M\}} \frac{|v(x,t)-v(y,t)|^p}{
|x-y|^{n+\sigma p}}\,dx\,dy
\end{eqnarray*}
and
\begin{eqnarray*}
&& C_\sharp \iint_{\R^{2n}\cap\{ |x-y|> M\}} \frac{|v(x,t)-v(y,t)|^p}{
|x-y|^{n+\tilde\sigma p}}\,dx\,dy\le C'
\iint_{\R^{2n}\cap\{ |x-y|> M\}} \frac{|v(x,t)|^p}{
|x-y|^{n+\tilde\sigma p}}\,dx\,dy\\
&&\qquad=
\frac{C''}{M^{\tilde\sigma p}}
\int_{\Omega} {|v(x,t)|^p}\,dx\le
\frac{C'''}{M^{\tilde\sigma p}}\,\|v\|^p_{L^{\tilde q}(\Omega)}(t)
\le \frac12\,\|v\|^p_{L^{\tilde q}(\Omega)}(t),
\end{eqnarray*}
as long as~$M$ is large enough. Here above, we have denoted by~$C'$, $C''$
and~$C'''$ suitable positive constants and used that~$\tilde q> p$
in order to use the H\"older's Inequality.
These inequalities and~\eqref{MAK0} imply that
\begin{equation}\label{MAK91} \frac12 \| v\|_{L^{\tilde q}(\Omega)}^{p}(t) \le C_\sharp\,
(1+M^{(\sigma-\tilde\sigma)p})
\iint_{\R^{2n}} \frac{|v(x,t)-v(y,t)|^p}{
|x-y|^{n+\tilde\sigma p}}\,dx\,dy.\end{equation}
We also have
that
\begin{equation}\label{MAK92}
\| v\|_{L^{q}(\Omega)}(t)\le C\,\| v\|_{L^{\tilde q}(\Omega)}(t),\end{equation}
for some~$C>0$,
in view of the H\"older's Inequality and the fact that~$\tilde q>q$.
Thanks to~\eqref{MAK91} and~\eqref{MAK92},
we have completed the proof of~\eqref{EM}.

Using~\eqref{DISV}, \eqref{EM} and the fact that~$u$ and~$v$ vanish outside~$\Omega$,
we see that, for every~$q\in[1,p_\sigma]$
when~$p\in\left(1,\frac{n}\sigma\right)$, and for every~$q\in[1,+\infty)$ when~$p
\in\left[\frac{n}\sigma,+\infty\right)$,
\begin{equation}\label{PIV}\begin{split}
&\left( \int_\Omega u^{(s-2+p)q/p}(x,t)\,dx\right)^{p/q}
\\=\,& \left( \int_{\R^n} v^{q}(x,t)\,dx\right)^{p/q}\\
=\,&
\| v\|_{L^{q}(\R^n)}^p(t)\\ \le\,& C_1\,
\iint_{\R^{2n}} \frac{|v(x,t)-v(y,t)|^p}{|x-y|^{n+\sigma p}}\,dx\,dy\\
\le\,&
C_2\,
\iint_{\R^{2n}} \frac{
|u(x,t)-u(y,t)|^{p-2}(u(x,t)-u(y,t))(u^{s-1}(x,t)-u^{s-1}(y,t))
}{|x-y|^{n+\sigma p}}\,dx\,dy\\
=\,&
2C_2\,
\iint_{\R^{2n}} \frac{
|u(x,t)-u(y,t)|^{p-2}(u(x,t)-u(y,t))\, u^{s-1}(x,t)
}{|x-y|^{n+\sigma p}}\,dx\,dy\\
=\,&
2C_2\,
\int_{\R^{n}} (-\Delta)^\sigma_p u(x,t)\, u^{s-1}(x,t)\,dx\\
=\,&
2C_2\,
\int_{\Omega} (-\Delta)^\sigma_p u(x,t)\, u^{s-1}(x,t)\,dx
,\end{split}\end{equation}
for some~$C_2>0$.

We also claim that when~$p\in\left(1,\frac{n}\sigma\right)$
and
\begin{equation}\label{stra}
s\geq\frac{n(2-p)}{\sigma p},
\end{equation}
it holds that
\begin{equation} \label{RANGE}\frac{sp}{s-2+p}\in [1,p_\sigma].\end{equation}
Indeed, $s-2+p>1-2+1=0$ and~$sp-s+2-p=s(p-1)+2-p>(p-1)+2-p=1$, which gives that~$
\frac{sp}{s-2+p}\ge1$. In addition,
$$ s(n-\sigma p)-n(s-2+p)=-s\sigma p+n(2-p)\le0,$$
thanks to~\eqref{stra}, which, recalling~\eqref{psigma}, says that~$\frac{sp}{s-2+p}\le
p_\sigma$. These considerations prove~\eqref{RANGE}.

{F}rom~\eqref{RANGE} it follows that if either~$s\geq\frac{n(2-p)}{\sigma p}$ or~$p\ge\frac{n}\sigma$
then we can choose~$q:=\frac{sp}{s-2+p}$ in~\eqref{PIV}. Consequently, we have that
\begin{equation}\label{AIJ}
\left( \int_\Omega u^{s}(x,t)\,dx\right)^{\frac{s-2+p}{s} }
\le 2C_2\,
\int_{\Omega} (-\Delta)^\sigma_p u(x,t)\, u^{s-1}(x,t)\,dx.\end{equation}
This says that~\eqref{SA} is satisfied with~$\gamma:=p-1$.
Hence, we are in position of exploiting~\eqref{v4}, obtaining that,
if either~$s\geq\frac{n(2-p)}{\sigma p}$ or~$p\ge\frac{n}\sigma$,
\begin{equation}\label{PERO} \| u\|_{L^s(\Omega)}(t)\le\frac{C_\star}{1+t^{\alpha/(p-1)}}.\end{equation}
We also observe that when~$s\in\left(1,\frac{n(2-p)}{\sigma p}\right)$
we have that
$$ \|u\|_{L^s(\Omega)}\le \hat C\,\|u\|_{L^{ \frac{n(2-p)}{\sigma p} }(\Omega)},$$
thanks to the H\"older's Inequality.
This and~\eqref{PERO} imply~\eqref{PERO2} for all~$s>1$ and~$p>1$.~\hfill$\Box$

\subsection{Proof of Theorem~\ref{THp2}}

The main idea is to use~\eqref{AIJ} for each 
index~$j\in\{ 1,\dots,N\}$.
That is, we fix
$$ \tilde s:=\max\left\{ s,\,\frac{n(2-p_1)}{\sigma_1 p_1},\dots,\frac{n(2-p_N)}{\sigma_N p_N} \right\}$$
and we exploit~\eqref{AIJ} to write that
\begin{equation}\label{00p01}
\| u\|_{L^{\tilde s}(\Omega)}^{\tilde s-2+p_j}
(t)=
\left(\int_\Omega u^{\tilde s}(x,t)\,dx
\right)^{\frac{\tilde s-2+p_j}{\tilde s}}
\le C\,\int_\Omega (-\Delta)^{\sigma_j}_{p_j} u(x,t)\,u^{\tilde s-1}(x,t)\,dx,\end{equation}
for some~$C>0$.

We also observe that
\begin{equation}\label{00p0}
\| u\|_{L^{\tilde s}(\Omega)}(t)\le \| u\|_{L^{\tilde s}(\Omega)}(0).\end{equation}
Indeed, we have that
$$ \big(u(x,t)-u(y,t)\big)\big(u^{\tilde s-1}(x,t)-u^{\tilde s-1}(y,t)\big)\ge0$$
and therefore
\begin{equation}\label{0wudyfg0efuihg}
\begin{split} 
&2\int_\Omega u^{\tilde s-1}(x,t)\,\NN[u](x,t)\,dx \\=\,&
2\sum_{j=1}^N \beta_j
\int_\Omega u^{\tilde s-1}(x,t)\,(-\Delta)_{p_j}^{\sigma_j} u(x,t)\,dx\\
=\,&2\sum_{j=1}^N \beta_j
\iint_{\R^{2n}} \frac{\big| u(x,t)-u(y,t)\big|^{p-2}
\big(u(x,t)-u(y,t)\big)}{|x-y|^{n+\sigma_j p_j}}
\,u^{\tilde s-1}(x,t)\,dx\,dy
\\=\,&
\sum_{j=1}^N \beta_j
\iint_{\R^{2n}} \frac{\big| u(x,t)-u(y,t)\big|^{p-2}
\big(u(x,t)-u(y,t)\big)\big(u^{\tilde s-1}(x,t)-u^{\tilde s-1}(y,t)\big)
}{|x-y|^{n+\sigma_j p_j}}\,dx\,dy
\\ \ge\,&0.\end{split}\end{equation}
Furthermore, from~\eqref{21a}, we know that
$$ \| u\|_{L^{\tilde s}(\Omega)}^{\tilde s-1}(t)\;\partial^\alpha_t 
\| u\|_{L^{\tilde s}(\Omega)}(t)\le -\int_\Omega u^{\tilde s-1}(x,t)\,\NN[ u](x,t)\,dx.$$
This and~\eqref{0wudyfg0efuihg}
give that
\begin{equation}\label{MAw9s:2} \| u\|_{L^{\tilde s}(\Omega)}^{\tilde s-1}(t)\;\partial^\alpha_t 
\| u\|_{L^{\tilde s}(\Omega)}(t)\le0.\end{equation}
We now observe that if~$\mu>0$, $f\ge0$, and
\begin{equation}\label{MAw9s}
{\mbox{$f^\mu(t)\,\partial^\alpha_t f(t)\le0$ with~$f(0)>0$, then~$\partial^\alpha_t f(t)\le0$.}}
\end{equation}
We prove\footnote{After this work was completed, arxived and submitted,
the very interesting preprint~\cite{DP} became available: with respect to this,
we mention that
formula~\eqref{MAw9s} here also follows from Lemma~2.1 in~\cite{DP}.}
this by contradiction, supposing that~$\partial^\alpha_t f(t_\star)>0$ for some~$t_\star>0$.
Hence, we find an open interval~$(a_\star,b_\star)$, with~$0<a_\star<t_\star$
such that~$\partial^\alpha_t f(t)>0$ for all~$t\in(a_\star,b_\star)$,
$f(t)=0$ for all~$t\in(a_\star,b_\star)$, and
$f(t)>0$ for all~$t\in[0,a_\star)$. This gives that $f(a_\star)=0$. Now, from~\eqref{K-DER}, integrating by parts twice we obtain that
\begin{eqnarray*}
\partial_t^{\alpha} f(t)&=&
\frac{1}{\alpha-1}\frac{d}{dt}
\int_0^t \big( f(\tau)-f(0)\big)\frac{d}{d\tau} (t-\tau)^{1-\alpha }\,d\tau
\\ &=&
\frac{1}{1-\alpha}\frac{d}{dt}\left[
\int_0^t \frac{d}{d\tau} \big( f(\tau)-f(0)\big)(t-\tau)^{1-\alpha }\,d\tau
\right]
\\ &=&
\frac{1}{1-\alpha}\frac{d}{dt}\left[
\int_0^t \dot f(\tau)(t-\tau)^{1-\alpha }\,d\tau
\right]
\\ &=&
\int_0^t \dot f(\tau)(t-\tau)^{-\alpha }\,d\tau
\\ &=&
\int_0^t 
\frac{d}{d\tau} \big( f(\tau)-f(t)\big)
(t-\tau)^{-\alpha }\,d\tau\\
&=&\frac{f(t)-f(0)}{t^\alpha}
+\alpha\int_0^t 
\frac{ f(t)-f(\tau)}{
(t-\tau)^{1+\alpha }}\,d\tau
\end{eqnarray*}
Consequently,
\begin{eqnarray*}0&\le&
\lim_{t\searrow a_\star} \partial_t^{\alpha} f(t)\\&=&
\partial_t^{\alpha} f(a_\star)\\&=&
\frac{f(a_\star)-f(0)}{a_\star^\alpha}
+\alpha\int_0^{a_\star} 
\frac{ f(a_\star)-f(\tau)}{
(a_\star-\tau)^{1+\alpha }}\,d\tau\\
&=&-\frac{f(0)}{a_\star^\alpha}
-\alpha\int_0^{a_\star} 
\frac{ f(\tau)}{
(a_\star-\tau)^{1+\alpha }}\,d\tau\\&<&0.
\end{eqnarray*}
This contradiction establishes~\eqref{MAw9s}.

By~\eqref{MAw9s:2}
and~\eqref{MAw9s}, we find that~$\partial^\alpha_t
\| u\|_{L^{\tilde s}(\Omega)}(t)\le0$.
{F}rom this,
we obtain~\eqref{00p0}
by inverting the Caputo derivative
by a Volterra
integral kernel
(see e.g. formula~(2.61) in~\cite{ABA};
alternatively, one could also
use the
comparison principle, e.g. Lemma~2.6 in~\cite{VER}).

Then, using~\eqref{00p01} and~\eqref{00p0}, we conclude that
$$ \| u\|_{L^{\tilde s}(\Omega)}^{\tilde s-2+p_{\max}}
(t)\le C'\,\| u\|_{L^{\tilde s}(\Omega)}^{\tilde s-2+p_j}\le
C''\,\int_\Omega (-\Delta)^{\sigma_j}_{p_j} u(x,t)\,u^{s-1}(x,t)\,dx,$$
for some~$C'$, $C''>0$ and so, multiplying by~$\beta_j>0$
and summing up over~$j\in\{1,\dots,N\}$,
$$ \| u\|_{L^{\tilde s}(\Omega)}^{\tilde s-2+p_{\max}}
(t)
\le C'''\,\int_\Omega \sum_{j=1}^N \beta_j\,(-\Delta)^{\sigma_j}_{p_j} u(x,t)\,u^{\tilde s-1}(x,t)\,dx.$$
This says that~\eqref{SA} is satisfied with~$s$ replaced 
by~$\tilde s$ and~$\gamma:=p_{\max}-1$.
Accordingly, we can exploit~\eqref{v4} and find that
\begin{equation}\label{ls67kd}
\| u\|_{L^{\tilde s}(\Omega)}(t)\le\frac{C_\star}{1+t^{\frac\alpha{p_{\max}-1}}}
.
\end{equation}
Since, by H\"older's Inequality
and the fact that~$s\le\tilde s$, we have that~$\| u\|_{L^{s}(\Omega)}\le C \| u\|_{L^{\tilde s}(\Omega)}$,
for some~$C>0$,
we deduce from~\eqref{ls67kd}
that~\eqref{PERO22} holds true, as desired.~\hfill$\Box$

\subsection{Proof of Theorem~\ref{THp3}}

We fix~$j\in\{1,\dots,n\}$ 
and~$(\rho_1,\dots,\rho_{j-1},\rho_{j+1},\dots,\rho_n)\in\R^{n-1}$
and denote~$\Omega_j(\rho_1,\dots,\rho_{j-1},\rho_{j+1},\dots,\rho_n):=\Omega\cap \R_j(\rho_1,\dots,\rho_{j-1},\rho_{j+1},\dots,\rho_n)$,
where
$$ \R_j(\rho_1,\dots,\rho_{j-1},\rho_{j+1},\dots,\rho_n):=\big\{ 
(\rho_1,\dots,\rho_{j-1},0,\rho_{j+1},\dots,\rho_n)+
r e_j,\;\, r\in\R\big\}.$$
The function~$\R\ni \rho_j\mapsto u(\rho_1 e_1+\dots+\rho_n e_n,t)$ is supported inside the closure
of the bounded set~$\Omega_j(\rho_1,\dots,\rho_{j-1},\rho_{j+1},\dots,\rho_n)$ and, using~\eqref{AIJ} with~$p:=2$, we get that
\begin{equation*}
\begin{split}&
\int_{\R} u^{s}(\rho_1 e_1+\dots+\rho_n e_n,t)\,d\rho_j
=
\int_{\Omega_j(\rho_1,\dots,\rho_{j-1},\rho_{j+1},\dots,\rho_n)} u^{s}(\rho_1 e_1+\dots+\rho_n e_n,t)\,d\rho_j\\&\qquad\qquad
\le C\,
\int_{\R} (-\partial^2_{x_j})^{\sigma_j} u(\rho_1 e_1+\dots+\rho_n e_n,t)\, u^{s-1}(
\rho_1 e_1+\dots+\rho_n e_n,t)\,d\rho_j,\end{split}\end{equation*}
for some~$C>0$. 

We now integrate such inequality
over the other coordinates~$(\rho_1,\dots,\rho_{j-1},\rho_{j+1},\dots,\rho_n)$
and we thus obtain that
\begin{equation*}\begin{split}
\int_{\Omega} u^{s}(x,t)\,dx\,&=\int_{\R^n} u^{s}(x,t)\,dx\\&=
\int_{\R^n} u^{s}(\rho_1 e_1+\dots+\rho_n e_n,t)\,d\rho
\\&\le C\,
\int_{\R^n} (-\partial^2_{x_j})^{\sigma_j} u(
\rho_1 e_1+\dots+\rho_n e_n
,t)\, u^{s-1}(
\rho_1 e_1+\dots+\rho_n e_n
,t)\,d\rho
\\&= C\,
\int_{\Omega} (-\partial^2_{x_j})^{\sigma_j} u(
x,t)\, u^{s-1}(x,t)\,dx.
\end{split}\end{equation*}
We multiply this inequality by~$\beta_j>0$ and we sum over~$j$, and we find that
$$ \int_{\Omega} u^{s}(x,t)\,dx\le C'
\, \int_{\Omega} (-\Delta_{\beta})^\sigma u(x,t)\, u^{s-1}(x,t)\,dx,$$
for some~$C'>0$.
Hence~\eqref{SA} holds true with~$\gamma:=1$
and then the desired result follows from~\eqref{v4}.~\hfill$\Box$

\subsection{Proof of Theorem~\ref{VVV}}

It is convenient to
define~$\tilde u:=u^m$, $\tilde s:=1+\frac{s-1}{m}$ and~$
\tilde v:=\tilde u^{\tilde s/2}$. Let also~$v:=u^{\frac{m+s-1}{2}}$.
We remark that
\begin{equation}\label{LAO1a49}
\tilde u^{\tilde s-1}=u^{s-1},\quad{\mbox{ and }}\quad
\tilde v=u^{m \tilde s/2}=u^{(m+s-1)/2}=v.
\end{equation}
We also exploit~\eqref{DISV} with~$p:=2$
to the functions~$\tilde u$ and~$\tilde v$,
with exponent~$\tilde s$.
In this way we have that
\begin{equation*}
|\tilde v(x,t)-\tilde v(y,t)|^2\le C_0 
(\tilde u(x,t)-\tilde u(y,t))(\tilde u^{\tilde s-1}(x,t)-\tilde u^{\tilde s-1}(y,t)),
\end{equation*}
for some~$C_0>0$. This estimate and~\eqref{LAO1a49} give that
\begin{equation}\label{porousDISV}
|v(x,t)-v(y,t)|^2\le C_0 
(u^m(x,t)- u^m(y,t))(u^{s-1}(x,t)-u^{s-1}(y,t)).
\end{equation}
Also, exploiting formula~\eqref{EM} with~$p:=2$, we have that
\begin{equation}\label{EMpor} 
\left(\int_\Omega u^{\frac{q(m+s-1)}{2}}(x,t)\,dx\right)^{2/q}
=\| v\|_{L^{q}(\Omega)}^{2}(t) \le C_1\,
\iint_{\R^{2n}} \frac{|v(x,t)-v(y,t)|^2}{
|x-y|^{n+2\sigma }}\,dx\,dy,\end{equation}
for some~$C_1>0$, for every~$q\in
\left[1,\frac{2n}{n-2\sigma }\right]$
when~$2\in\left(1,\frac{n}\sigma\right)$, 
and for every~$q\in[1,+\infty)$ when~$2
\in\left[\frac{n}\sigma,+\infty\right)$.

Now we observe that when
\begin{equation}\label{67a}
s\ge m-1,\end{equation} it holds that
\begin{equation}\label{67}
\frac{2s}{m+s-1}\in\left[1,\frac{2n}{n-2\sigma }\right].
\end{equation}
Indeed, we have that
$$ s(n-2\sigma)-n(m+s-1)=
-2\sigma s -n(m-1)\le0,$$
which says that~$\frac{2s}{m+s-1}\leq\frac{2n}{n-2\sigma }$.
On the other hand, from~\eqref{67a}, we see that
$$ 2s-(m+s-1)=s-m+1\ge0,$$
giving that~$\frac{2s}{m+s-1}\ge1$.
This proves~\eqref{67}.

Now, by~\eqref{67}, when either~$2
\in\left[\frac{n}\sigma,+\infty\right)$ or~$s\ge m-1$,
we are allowed to
choose~$q:=\frac{2s}{m+s-1}$ in~\eqref{EMpor}, and conclude 
that
$$ \left(\int_\Omega u^{s}(x,t)\,dx\right)^{(m+s-1)/s}
\le C_1\,
\iint_{\R^{2n}} \frac{|v(x,t)-v(y,t)|^2}{
|x-y|^{n+2\sigma }}\,dx\,dy.$$
Combining this with~\eqref{porousDISV}, we find that
\begin{eqnarray*} 
\left(\int_\Omega u^{s}(x,t)\,dx\right)^{(m+s-1)/s}
&\le& C_2\,\iint_{\R^{2n}} \frac{
(u^m(x,t)- u^m(y,t))(u^{s-1}(x,t)-u^{s-1}(y,t))}{
|x-y|^{n+2\sigma}}\,dx\,dy\\
&=&
2C_2\,\iint_{\R^{2n}} \frac{
(u^m(x,t)- u^m(y,t))\,u^{s-1}(x,t)}{|x-y|^{n+2\sigma}}\,dx\,dy,\end{eqnarray*}
provided that either~$2
\in\left[\frac{n}\sigma,+\infty\right)$ or~$s\ge m-1$.

This says that, under these circumstances,
condition~\eqref{SA} is fulfilled with~$\gamma:=m$.
Therefore, we can exploit~\eqref{v4}
and obtain~\eqref{sdfbgssh82394:2}, 
provided that either~$2
\in\left[\frac{n}\sigma,+\infty\right)$ or~$s\ge m-1$.

Then, when~$2<\frac{n}\sigma$, we first establish~\eqref{sdfbgssh82394:2}
for a large exponent of the Lebesgue norm, and then we reduce it
by using the H\"older's Inequality.

This completes the proof of~\eqref{sdfbgssh82394:2}
in all the cases,
as desired.~\hfill$\Box$

\subsection{Proof of Theorem~\ref{FF}}

We claim that
\begin{equation}\label{DAF3S}
\begin{split}&F\left( \frac{u(x,t)-u(y,t)}{|x-y|}\right)\,\big(u^{s-1}(x,t)-u^{s-1}(y,t)\big)\\ \ge\;& c_0\,
\frac{\big( u(x,t)-u(y,t) \big)\big(u^{s-1}(x,t)-u^{s-1}(y,t)\big)}{|x-y|}\end{split}
\end{equation} for some~$c_0>0$.
To check this, we
observe that, by~\eqref{F}, $F$ is odd, hence we can reduce to the case in which
\begin{equation}\label{0wo20}u(x,t)\ge u(y,t).\end{equation}
Also, by~\eqref{F}, we see that when~$|r|$ is bounded, then~$F(r)\simeq r$ and therefore
by~\eqref{GBOU3} and~\eqref{0wo20},
$$ F\left( \frac{u(x,t)-u(y,t)}{|x-y|}\right)\ge c_0\,
\frac{u(x,t)-u(y,t)}{|x-y|},$$
and this implies \eqref{DAF3S}.

Now, we let~$v:=u^{s/2}$.
We use~\eqref{DISV} (with~$p:=2$)
and~\eqref{DAF3S} to deduce that
\begin{equation}\label{231}
\begin{split}&
2\iint_{\R^{2n}} \frac{1}{|y|^{n+\sigma}}
\,F\left( \frac{u(x,t)-u(x+y,t)}{|y|}\right)\, u^{s-1}(x,t)\,dx\,dy\\
=\;&
2\iint_{\R^{2n}} \frac{1}{|x-y|^{n+\sigma}}
\,F\left( \frac{u(x,t)-u(y,t)}{|x-y|}\right)\, u^{s-1}(x,t)\,dx\,dy
\\=\;&
\iint_{\R^{2n}} \frac{1}{|x-y|^{n+\sigma}}
\,F\left( \frac{u(x,t)-u(y,t)}{|x-y|}\right)\, u^{s-1}(x,t)\,dx\,dy
\\&\qquad+
\iint_{\R^{2n}} \frac{1}{|x-y|^{n+\sigma}}
\,F\left( \frac{u(y,t)-u(x,t)}{|x-y|}\right)\, u^{s-1}(y,t)\,dx\,dy
\\=\;&
\iint_{\R^{2n}} \frac{1}{|x-y|^{n+\sigma}}
\,F\left( \frac{u(x,t)-u(y,t)}{|x-y|}\right)\, \big(u^{s-1}(x,t)-u^{s-1}(y,t)\big)\,dx\,dy
\\ \ge\;&c_0\,\iint_{\R^{2n}}
\frac{\big( u(x,t)-u(y,t) \big)\big(u^{s-1}(x,t)-u^{s-1}(y,t)\big)}{|x-y|^{n+\sigma+1}}\,dx\,dy
\\ \ge\;&\frac{c_0}{C_0}\,\iint_{\R^{2n}}
\frac{\big| v(x,t)-v(y,t) \big|^2}{|x-y|^{n+2\sigma'}}\,dx\,dy,
\end{split}\end{equation}
where~$\sigma':=\frac{\sigma+1}2\in(0,1)$.

Furthermore, we recall~\eqref{EM}, used here with fractional exponent~$\sigma'$ and with~$p:=2$,
and we see that
\begin{equation}\label{FR}
\left(
\int_\Omega u^{sq/2}(x,t)\,dx\right)^{2/q}=
\| v\|_{L^{q}(\Omega)}^{2}(t) \le C_1\,
\iint_{\R^{2n}} \frac{\big|v(x,t)-v(y,t)\big|^2}{
|x-y|^{n+2\sigma'}}\,dx\,dy,\end{equation}
for some~$C_1>0$, for every~$q\in\left[1,\frac{2n}{n-2\sigma'}\right]$
when~$2\in\left(1,\frac{n}{\sigma'}\right)$, and for every~$q\in[1,+\infty)$ when~$2
\in\left[\frac{n}{\sigma'},+\infty\right)$. 

In any case, since~$\frac{2n}{n-2\sigma'}>2$,
we can always choose~$q:=2$ in~\eqref{FR} and conclude that
\begin{equation*}
\int_\Omega u^{s}(x,t)\,dx \le C_1\,
\iint_{\R^{2n}} \frac{\big|v(x,t)-v(y,t)\big|^2}{
|x-y|^{n+2\sigma'}}\,dx\,dy.\end{equation*}
Combining this with~\eqref{231}, we infer that
\begin{eqnarray*}\int_\Omega u^{s}(x,t)\,dx\le C_2\,\iint_{\R^{2n}} \frac{1}{|y|^{n+\sigma}}
\,F\left( \frac{u(x,t)-u(x+y,t)}{|y|}\right)\, u^{s-1}(x,t)\,dx\,dy\end{eqnarray*}
for some~$C_2>0$, which, together with~\eqref{CH6}, establishes~\eqref{SA} with~$\gamma:=1$.
This and~\eqref{v4} yield the thesis of Theorem~\ref{FF}.~\hfill$\Box$

\end{document}